\documentclass[11pt,twoside]{article}
\usepackage{multicol}
\usepackage[greek,english]{babel}
\usepackage[vcentering,dvips]{geometry}
\usepackage{amsfonts}
\usepackage{amssymb}
\usepackage{eucal}
\usepackage{amsxtra}
\usepackage{color}
\usepackage{setspace}
\usepackage{times}
\usepackage{titlesec}
\usepackage{sectsty}
\usepackage{setspace}
\usepackage{fancyhdr}
\usepackage{fancyheadings}
\allsectionsfont{\large}
\usepackage{url}
\usepackage{doi}
\usepackage{hyperref}
\hypersetup{
    colorlinks,%
    citecolor=black,%
    filecolor=black,%
    linkcolor=black,%
    urlcolor=black
}

\addtocounter{section}{0} \numberwithin{equation}{section}
\newtheorem{theorem}{Theorem}[section]
\newtheorem{proposition}[theorem]{Proposition}

\newtheorem{lemma}[theorem]{Lemma}

\hypersetup{urlcolor=blue }

\begin{document}
\pagestyle{fancy}
\fancyhead[LE,RO]{\thepage}
\renewcommand{\sectionmark}[1]{\markright{\emph{ \thesection.\ #1}}{}}
\fancyhead[RE]{\leftmark}
\fancyhead[LO]{\rightmark}
\cfoot{}

\title{\textbf{Parallel *-Ricci tensor of real hypersurfaces in $\mathbb{C}P^{2}$ and $\mathbb{C}H^{2}$}}

\author{\textsc{George Kaimakamis and Konstantina Panagiotidou}}
\date{}
\maketitle

{\hspace{-15pt}\textsc{Abstract:}\small   In this paper the idea of studying real hypersurfaces in non-flat complex space forms, whose *-Ricci tensor satisfies geometric conditions is presented. More precisely, the non-existence of three dimensional real hypersurfaces in non-flat complex space forms with parallel *-Ricci tensor is proved. At the end of the paper ideas for further research on $^{*}$-Ricci tensor are given.}
\begin{flushleft}
\small{\emph{Keywords}: Real hypersurface, Parallel, *-Ricci tensor, Complex projective plane, Complex hyperbolic plane.\\}
\end{flushleft}
\begin{flushleft}
\small{\emph{Mathematics Subject Classification }(2010):  Primary 53B20; Secondary 53C15, 53C25.}
\end{flushleft}

\rhead[\centering{G. Kaimakamis and K. Panagiotidou}]{\thepage}
\lhead[\thepage]{\centering{Parallel *-Ricci tensor}}

\section{Introduction}
A \emph{complex space form} is an $n$-dimensional Kaehler manifold	of constant holomorphic sectional curvature \emph{c}. A complete and simply connected complex space form is complex analytically isometric to
\begin{itemize}
  \item  a complex projective space $\mathbb{C}P^{n}$ if $c>0$,
  \item  a complex Euclidean space $\mathbb{C}^{n}$ if $c=0$,
  \item  or a complex hyperbolic space $\mathbb{C}H^{n}$ if $c<0$.
\end{itemize}
The symbol $M_{n}(c)$ is used to denote the complex projective space $\mathbb{C}P^{n}$  and complex hyperbolic space $\mathbb{C}^{n}$, when it is not necessary to distinguish them. Furthermore, since $c\neq0$ in previous two cases the notion of non-flat complex space form refers to both them.

Let \emph{M} be a real hypersurface in a non-flat complex space form. An almost contact metric structure $(\varphi, \xi, \eta, g)$ is defined on \emph{M} induced from the Kaehler metric \emph{G} and the complex structure \emph{J} on $M_{n}(c)$. The \emph{structure vector field} $\xi$ is called \emph{principal} if $A\xi=\alpha\xi$, where \emph{A} is the shape operator of \emph{M} and $\alpha=\eta(A\xi)$ is a smooth function. A real hypersurface is called \emph{Hopf hypersurface}, if $\xi$ is principal and $\alpha$ is called \emph{Hopf principal curvature}.

The \emph{Ricci tensor} $S$ of a Riemannian manifold is a tensor field of type (1,1) and is given by
\[g(SX,Y)=trace\{Z\mapsto R(Z,X)Y\}.\]
If the Ricci tensor of a Riemannian manifold satisfies the relation
\[S=\lambda g,\]
where $\lambda$ is a constant, then it is called \emph{Einstein}.

Real hypersurfaces in non-flat complex space forms have been studied in terms of their Ricci tensor $S$, when it satisfies certain geometric conditions extensively. Different types of parallelism or invariance of the Ricci tensor are issues of great importance in the study of real hypersurfaces.

In \cite{Ki1989} it was proved the non-existence of real hypersurfaces in non-flat complex space forms $M_{n}(c)$, $n\geq3$ with parallel Ricci tensor, i.e. $(\nabla_{X}S)Y=0$, for any $X$, $Y$ $\in$ $TM$. In \cite{Kim2004} Kim extended the result of non-existence of real hypersurfaces with parallel Ricci tensor in case of three dimensional real hypersurfaces. Another type of parallelism which was studied is that of $\xi$-parallel Ricci tensor, i.e. $(\nabla_{\xi}S)Y=0$ for any $Y$ $\in$ $TM$. More precisely in \cite{KimMaed1991} Hopf hypersurfaces in non-flat complex space forms with constant mean curvature and $\xi$-parallel Ricci tensor were classified. More details on the study of Ricci tensor of real hypersurfaces are included in Section 6 of \cite{NR1}.

Motivated by Tachibana, who in \cite{Tach} introduced the notion of \emph{*-Ricci tensor} on almost Hermitian manifolds, in \cite{Ham} Hamada defined the \emph{*-Ricci tensor} of real hypersurfaces in non-flat complex space forms by
\[g(S^{*}X,Y)=\frac{1}{2}(trace\{\varphi \circ R(X,\varphi Y)\}), \;\;\mbox{ for $X$, $Y$ $\in$ $TM$}.\]
The $^{*}$-Ricci tensor $S^{*}$ is a tensor field of type (1,1) defined on real hypersurfaces. Taking into account the work that so far has been done in the area of studying real hypersurfaces in non-flat complex space forms in terms of their tensor fields, the following issue raises naturally:

\emph{The study of real hypersurfaces in terms of their $^{*}$-Ricci tensor $S^{*}$, when it satisfies certain geometric conditions.}

In this paper three dimensional real hypersurfaces in $\mathbb{C}P^{2}$ and $\mathbb{C}H^{2}$ equipped with parallel $^{*}$-Ricci tensor are studied. Therefore, the following condition is satisfied
\begin{eqnarray}\label{BS}
(\nabla_{X}S^{*})Y=0,\;\;\mbox{$X$, $Y$ $\in$ $TM$.}
\end{eqnarray}

More precisely the following Theorem is proved.

\begin{prob}\label{MT}
There do not exist real hypersurfaces in $\mathbb{C}P^{2}$ and in $\mathbb{C}H^{2}$, whose *-Ricci tensor is parallel.
\end{prob}

The paper is organized as follows: In Section 2 preliminaries relations for real hypersurfaces in non-flat complex space forms are presented. In Section 3 the proof of Main Theorem is provided. Finally, in Section 4 ideas for further research on the above issue are included.

\section{Preliminaries}
Throughout this paper all manifolds, vector fields etc are assumed to be of class $C^{\infty}$ and all manifolds are assumed to be connected. Furthermore, the real hypersurfaces $M$ are supposed to be without boundary.

Let $M$ be a real hypersurface immersed in a non-flat complex space form $(M_{n}(c),G)$ with complex structure $J$ of constant holomorphic sectional curvature $c$. Let $N$ be a locally defined unit normal vector field on $M$ and $\xi=-JN$ the structure vector field of $M$.

For a vector field $X$ tangent to $M$ the following relation holds
\[JX=\varphi X+\eta(X)N,\]
where $\varphi X$ and $\eta(X)N$ are the tangential and the normal
component of $JX$ respectively. The Riemannian connections
$\overline{\nabla}$ in $M_{n}(c)$ and $\nabla$ in $M$ are related
for any vector fields $X$, $Y$ on $M$ by
\[\overline{\nabla}_{X}Y=\nabla_{X}Y+g(AX,Y)N,\]
where $g$ is the Riemannian metric induced from the metric $G$.

The shape operator $A$ of the real hypersurface $M$ in $M_{n}(c)$ with respect to $N$ is given by
\[\overline{\nabla}_{X}N=-AX.\]
The real hypersurface $M$ has an almost contact metric structure $(\varphi,\xi,\eta, g)$ induced from the complex structure $J$ on $M_{n}(c)$, where $\varphi$  is the \emph{structure tensor} and it is a tensor field of type (1,1). Moreover, $\eta$ is an
1-form on $M$ such that
\[g(\varphi X,Y)=G(JX,Y),\hspace{20pt}\eta(X)=g(X,\xi)=G(JX,N).\]
Furthermore, the following relations hold
\begin{eqnarray}\label{eq-1}
\varphi^{2}X=-X+\eta(X)\xi,\hspace{20pt}
\eta\circ\varphi=0,\hspace{20pt} \varphi\xi=0,\hspace{20pt}
\eta(\xi)=1,\nonumber\\
\hspace{20pt}
g(\varphi X,\varphi
Y)=g(X,Y)-\eta(X)\eta(Y),\hspace{10pt}g(X,\varphi Y)=-g(\varphi
X,Y).\label{eq-2}\nonumber\
\end{eqnarray}
Since $J$ is complex structure implies $\nabla J=0$. The last relation leads to
\begin{eqnarray}\label{eq-3}
\nabla_{X}\xi=\varphi
AX,\hspace{20pt}(\nabla_{X}\varphi)Y=\eta(Y)AX-g(AX,Y)\xi.
\end{eqnarray}
    The ambient space $M_{n}(c)$ is of constant holomorphic sectional
curvature $c$ and this results in the Gauss and Codazzi equations to be given respectively by
\begin{eqnarray}\label{eq-4}
R(X,Y)Z=\frac{c}{4}[g(Y,Z)X-g(X,Z)Y+g(\varphi Y ,Z)\varphi
X
\end{eqnarray}
$$-g(\varphi X,Z)\varphi Y-2g(\varphi X,Y)\varphi
Z]+g(AY,Z)AX-g(AX,Z)AY,$$
\begin{eqnarray}\label{eq-5}
\hspace{10pt}
(\nabla_{X}A)Y-(\nabla_{Y}A)X=\frac{c}{4}[\eta(X)\varphi
Y-\eta(Y)\varphi X-2g(\varphi X,Y)\xi]\nonumber\,
\end{eqnarray}
where $R$ denotes the Riemannian curvature tensor on $M$ and $X$, $Y$, $Z$ are any vector fields on $M$.

The tangent space $T_{P}M$ at every point $P$ $\in$ $M$  can be decomposed as
\begin{eqnarray}
T_{P}M=span\{\xi\}\oplus \mathbb{D},\nonumber\
\end{eqnarray}
where $\mathbb{D}=\ker\eta=\{X\;\in\;T_{P}M:\eta(X)=0\}$ and is called \emph{holomorphic distribution}.
 Due to the above decomposition the vector field $A\xi$ can be written
 \begin{eqnarray}\label{eq-7}
 A\xi=\alpha\xi+\beta U,
 \end{eqnarray}
 where $\beta=|\varphi\nabla_{\xi}\xi|$ and
 $U=-\frac{1}{\beta}\varphi\nabla_{\xi}\xi\;\in\;\ker(\eta)$ provided
 that $\beta\neq0$.
 
Since the ambient space $M_{n}(c)$ is of constant holomorphic sectional curvature $c$ following similar calculations to those in Theorem 2  in \cite{IvRy} and taking into account relation (\ref{eq-4}), it is proved that the \emph{*-Ricci tensor $S^{*}$} of $M$ is given by
\begin{eqnarray}\label{eq-8}
S^{*}= - [\frac{cn}{2}\varphi^{2}+(\varphi A)^{2}].
\end{eqnarray}

\section{Proof of Main Theorem}

Let \emph{M} be a non-Hopf hypersurface in $\mathbb{C}P^{2}$ or $\mathbb{C}H^{2}$, i.e. $M_{2}(c)$. Then the following relations hold on every non-Hopf three-dimensional real hypersurface in $M_{2}(c)$.

\begin{lemma}\label{lemma-1}
Let M be a real hypersurface in $M_{2}(c)$. Then the following relations hold on M
\begin{eqnarray}
&&AU=\gamma U+\delta\varphi U+\beta\xi,\;\;A\varphi U=\delta U+\mu\varphi U,\label{B1}\\
&&\nabla_{U}\xi=-\delta U+\gamma\varphi U,\;\;\nabla_{\varphi U}\xi=-\mu U+\delta\varphi U,\;\;\nabla_{\xi}\xi=\beta\varphi U,\label{B2}\\
&&\nabla_{U}U=\kappa_{1}\varphi U+\delta\xi,\;\;\nabla_{\varphi U}U=\kappa_{2}\varphi U+\mu\xi,\;\;\nabla_{\xi}U=\kappa_{3}\varphi U,\label{B3}\\
&&\nabla_{U}\varphi U=-\kappa_{1}U-\gamma\xi,\;\;\nabla_{\varphi U}\varphi U=-\kappa_{2}U-\delta\xi,\;\;\nabla_{\xi}\varphi U=-\kappa_{3}U-\beta\xi,\label{B4}
\end{eqnarray}
where $\gamma,$ $\delta,$ $\mu$, $\kappa_{1}$, $\kappa_{2}$, $\kappa_{3}$ are smooth functions on M and $\{U,\varphi U,\xi\}$ is an orthonormal basis of M.
\end{lemma}
For the proof of the above Lemma see \cite{PX3}

Let $M$ be a real hypersurface in $M_{2}(c)$, i.e. $\mathbb{C}P^{2}$ or $\mathbb{C}H^{2}$, whose $^{*}$-Ricci tensor satisfies relation (\ref{BS}), which is more analytically written
\begin{eqnarray}\label{BS-1}
\nabla_{X}(S^{*}Y)=S^{*}(\nabla_{X}Y),\;\;\mbox{$X$, $Y$ $\in$ $TM$.}
\end{eqnarray}
We consider the open subset $\mathcal{N}$ of \emph{M} such that
\begin{align*}
\mathcal{N}=\{P\;\;\in\;\;M:\beta\neq0,\;\mbox{in a neighborhood of $P$}\}.
\end{align*}
In what follows we work on the open subset $\mathcal{N}$.

On $\mathcal{N}$ relation (\ref{eq-7}) and relations (\ref{B1})-(\ref{B4}) of Lemma \ref{lemma-1} hold. So relation (\ref{eq-8}) for $X$ $\in$ $\{U, \varphi U,\xi\}$ taking into account $n=2$  and relations (\ref{eq-7}) and (\ref{B1}) yields
\begin{eqnarray}\label{C1}
S^{*}\xi=\beta\mu U-\beta\delta\varphi U,\;\;\;\;S^{*}U=(c+\gamma\mu-\delta^{2})U\;\;\mbox{and}\;\;S^{*}\varphi U=(c+\gamma\mu-\delta^{2})\varphi U.
\end{eqnarray}
The inner product of relation (\ref{BS-1}) for $X=Y=\xi$ with $\xi$ due to the first and the third of (\ref{C1}), the first of (\ref{eq-3}) for $X=\xi$ and the third of relations (\ref{B3}) and (\ref{B4}) implies
\begin{eqnarray}\label{C2}
\delta=0.
\end{eqnarray}
Moreover, the inner product of relation (\ref{BS-1}) for $X=\varphi U$ and $Y=\xi$ with $\xi$ because of (\ref{C2}), the first of (\ref{eq-3}) for $X=\varphi U$, the first and the second of (\ref{C1}) and the second of (\ref{B3}) results in
\begin{eqnarray}\label{C3}
\mu=0.\nonumber\
\end{eqnarray}
Finally, the inner product of relation (\ref{BS-1}) for $X=\xi$ and $Y=\varphi U$ with $\xi$ taking into account $\mu=\delta=0$, the first and the third of (\ref{C1}) and the last relation of (\ref{B4}) leads to
\[c=0,\]
which is a contradiction. So the open subset $\mathcal{N}$ is empty and we lead to the following Proposition.

\begin{proposition}
Every real hypersurface in $M_{2}(c)$ whose $^{*}$-Ricci tensor is parallel, is a Hopf hypersurface.
\end{proposition}

Since $M$ is a Hopf hypersurface, the structure vector field $\xi$ is an eigenvector of the shape operator, i.e. $A\xi=\alpha\xi$. Due to Theorem 2.1 in \cite{NR1} $\alpha$ is constant. We consider a point $P$ $\in$ \emph{M} and choose a unit principal vector field \emph{W} $\in$ $\mathbb{D}$ at \emph{P}, such that $AW=\lambda W$ and $A\varphi W=\nu\varphi W$. Then $\{W, \varphi W, \xi\}$ is a local orthonormal basis and the following relation holds (Corollary 2.3 \cite{NR1})
\begin{eqnarray}\label{Basic-Relation}
\lambda\nu=\frac{\alpha}{2}(\lambda+\nu)+\frac{c}{4}.
\end{eqnarray}

The first of relation (\ref{eq-3}) and relation (\ref{eq-8}) for $X$ $\in$ $\{W, \varphi W, \xi\}$ because of $A\xi=\alpha\xi$, $AW=\lambda W$ and $A\varphi W=\nu\varphi W$ implies respectively
\begin{eqnarray}
\nabla_{W}\xi=\lambda\varphi W\;\;\mbox{and}\;\;\nabla_{\varphi W}\xi=-\nu W\label{D1}
\end{eqnarray}
\begin{eqnarray}
S^{*}\xi=0,\;\;\;\;S^{*}W=(c+\lambda\nu)W\;\;\mbox{and}\;\;S^{*}\varphi W=(c+\lambda\nu)\varphi W\label{D2}.
\end{eqnarray}
Relation (\ref{BS-1}) for $X=W$ and $Y=\xi$ because of the first of (\ref{D1}) and the first and third relation of ({\ref{D2}) yields
\[\lambda(c+\lambda\nu)=0.\]

Suppose that $(c+\lambda\nu)\neq0$ then the above relation results in $\lambda=0$. Moreover, relation (\ref{BS-1}) for $X=\varphi W$ and $Y=\xi$ because of the second of (\ref{D1}) and the first and second relation of ({\ref{D2}) yields
\[\nu=0.\]
Substitution of $\lambda=\nu=0$ in (\ref{Basic-Relation}) results in $c=0$, which is a contradiction. So relation $c=-\lambda\nu$ holds. The last one implies $\lambda\nu\neq0$ since $c\neq0$.

Suppose that $\lambda=-\frac{c}{\nu}$. Substitution of the last one in (\ref{Basic-Relation}) leads to
\begin{eqnarray}\label{equation}
2\alpha\nu^{2}+5c\nu-2\alpha c=0.
\end{eqnarray}

In case of $\mathbb{C}P^{2}$ we have that $c=4$ and from equation (\ref{equation}) there is always a solution for $\nu$. So 
$\nu$ is constant and $\lambda$ will be also constant. Therefore, the real hypersurface has three distinct constant eigenvalues. The latter results in $M$ being a real hypersurface of type ($B$), i.e. a tube of radius $r$ over complex quadric. Substitution of the eigenvalues of type ($B$) in $\lambda\nu=-c$ leads to a contradiction. So no real hypersurface in $\mathbb{C}P^{2}$ has parallel $^{*}$-Ricci tensor (eigenvalues can be found in \cite{NR1}).

In case of $\mathbb{C}H^{2}$ we have that $c=-4$ and from equation (\ref{equation}) there is a solution for $\nu$ if $0\leq\alpha^{2}\leq\frac{25}{4}$. If $\alpha=0$ equation (\ref{equation}) implies $c\nu=0$, which is impossible. So there is a solution for $\nu$ if  $0<\alpha^{2}\leq\frac{25}{4}$ and $\nu$ will be constant. The latter results in that $\lambda$ is also constant and so the real hypersurface is of type ($B$), i.e. a tube of radius $r$ around totally geodesic $\mathbb{R}H^{n}$. Substitution of the eigenvalues of type ($B$) in $\lambda\nu=-c$ leads to a contradiction and this completes the proof of our Main Theorem  (eigenvalues can be found in \cite{Ber}).

\section{Discussion-Open Problems}
In this paper three dimensional real hypersurfaces in non-flat complex space forms with parallel $^{*}$-Ricci tensor are studied and the non-existence of them is proved. Therefore, a question which raises in a natural way is

\emph{Are there real hypersurfaces in non-flat complex space forms of dimension greater than three with parallel $^{*}$-Ricci tensor?}

Generally, the next step in the study of real hypersurfaces in non-flat complex space forms is to study them when a tensor field $P$ type (1,1) of them satisfies other types of parallelism such as the $\mathbb{D}$-parallelism or $\xi$-parallelism. The first one implies that $P$ is parallel in the direction of any  vector field $X$ orthogonal to $\xi$, i.e. $(\nabla_{X}P)Y=0$, for any $X$ $\in$ $\mathbb{D}$, and the second one implies that $P$ is parallel in the direction of the structure vector $\xi$, i.e. $(\nabla_{\xi}P)Y=0$. So the questions which should be answered are the following

\emph{Are there real hypersurfaces in non-flat complex space forms whose $^{*}$-Ricci tensor satisfies the condition of $\mathbb{D}$-parallelism or $\xi$-parallelism?}

Finally, other types of parallelism play important role in the study of real hypersurfaces is that of semi-parallelism and pseudo-parallelism. A tensor field $P$ of type (1, \emph{s}) is said to be
\emph{semi-parallel} if it satsfies $R\cdot P=0$, where $R$ is the Riemannian curvature tensor and acts as a derivation on $P$. Moreover, $P$ is said to be
\emph{pseudo-parallel} if there exists a function \emph{L} such that $R(X,Y)\cdot P=L \{(X\wedge Y)\cdot P\},$ where $(X\wedge Y)Z=g ( Y, Z )X-g( Z, X)Y$. So the questions are:

\emph{Are there real hypersurfaces in non-flat complex space forms with semi-parallel or pseudo-parallel $^{*}$-Ricci tensor?}

The importance of answering the above question lays in the fact that the class of real hypersurfaces with parallel $*$-Ricci tensor is included in the class of real hypersurfaces with semi-parallel $^{*}$-Ricci tensor. Furthermore, the last one is included in the class of real hypersurfaces with pseudo-parallel $^{*}$-Ricci tensor.

\end{document}